\documentclass[12pt,amssymb,openbib]{article}
\usepackage{amssymb}
\usepackage{epsf}
\usepackage{graphicx}

\usepackage{latexsym,epsf}
\title{The Generalized-Euler-Constant Function $\gamma(z)$ and \\ a Generalization of  Somos's Quadratic Recurrence 
Constant}
\author{Jonathan Sondow \\ 209 West 97th Street \\ New York, NY 10025, USA \\ {\tt jsondow@alumni.princeton.edu} \\
\\ Petros Hadjicostas\\ Department of Mathematics and Statistics \\ Texas Tech University \\ Box 41042 \\ Lubbock, TX 
79409-1042, USA \\ {\tt petros.hadjicostas@ttu.edu}}
\date{\today}
\newcommand{\R}{{\mathbb R}}
\newcommand{\Z}{{\mathbb Z}}

\newcommand{\Arg}{{\rm Arg}\,}

\newcommand{\be}{\begin{equation}}
\newcommand{\ee}{\end{equation}}
\newtheorem{lemma}{Lemma}
\newtheorem{theorem}{Theorem}

\newtheorem{example}{Example}
\newtheorem{remark}{Remark}

\newtheorem{definition}{Definition}
\newtheorem{corollary}{Corollary}

\newcommand{\bl}{\begin{lemma}}
\newcommand{\el}{\end{lemma}}
\newcommand{\bt}{\begin{theorem}}
\newcommand{\et}{\end{theorem}}
\newcommand{\bcor}{\begin{corollary}}
\newcommand{\ecor}{\end{corollary}}
\newcommand{\C}{{\mathbb{C}}}

\begin{document}
\maketitle

\begin{abstract}
We define the generalized-Euler-constant function $\gamma(z)=\sum_{n=1}^{\infty} z^{n-1}
\left(\frac{1}{n}-\log \frac{n+1}{n}\right)$ when $|z|\leq 1$. Its values include both Euler's constant
$\gamma=\gamma(1)$ and the ``alternating Euler constant" $\log\frac{4}{\pi}=\gamma(-1)$.
We extend Euler's two zeta-function series for $\gamma$ to polylogarithm series for $\gamma(z)$.
Integrals for $\gamma(z)$ provide its analytic continuation to $\C-[1,\infty)$. We prove several other
formulas for $\gamma(z)$, including two functional equations; one is an inversion relation between $\gamma(z)$ and 
$\gamma(1/z)$.  We generalize Somos's quadratic recurrence
constant and sequence to cubic and other degrees, 
give asymptotic estimates, and show relations to $\gamma(z)$ and to an infinite nested
radical due to Ramanujan. We calculate $\gamma(z)$ and $\gamma'(z)$ at roots of unity; in particular, $\gamma'(-1)$
involves the Glaisher-Kinkelin constant $A$. Several related series, infinite products, and double integrals are
evaluated. The methods used involve the Kinkelin-Bendersky hyperfactorial $K$ function, the Weierstrass products
for the gamma and Barnes $G$ functions, and Jonqui\`{e}re's relation for the polylogarithm. 
\end{abstract}

\vspace{.1in}

\begin{center}
{\bf Contents}
$$\begin{array}{llr}
1. & \mbox{Introduction} & \pageref{zerosection} \\
2. & \mbox{The Generalized-Euler-Constant Function $\gamma(z)$} & \pageref{twosection} \\
3. & \mbox{A Generalization of Somos's Quadratic Recurrence Constant} & \pageref{threesection} \\
4. & \mbox{Calculation of $\gamma(z)$ at Roots of Unity} & \pageref{onesection}\\
5. & \mbox{The Hyperfactorial $K$ Function and the Derivative of $\gamma(z)$} & \pageref{foursection} \\
 & \mbox{References} & \pageref{refref}
\end{array}$$
\end{center}

\noindent {\large {\bf 1. Introduction}}\label{zerosection}

\vspace{.2in}

In this paper we introduce and study the generalized-Euler-constant function $\gamma(z)$, defined by 
$$\gamma(z)=\sum_{n=1}^{\infty} z^{n-1} \left(\frac{1}{n}-\log\frac{n+1}{n}\right)=
\int_0^1\int_0^1 \frac{1-x}{(1-xyz)(-\log xy)}\,dx\,dy,$$
where the series converges when $|z|\leq 1$, and the integral gives the analytic continuation for $z \in \C -[1,\infty).$
The function $\gamma(z)$ generalizes both  Euler's constant $\gamma=\gamma(1)$ and the ``alternating Euler constant"
$\log \frac{4}{\pi}=\gamma(-1)$ \cite{sonqqq}, \cite{quisond22}, where $\gamma$ is defined by the limit
\begin{equation}\label{defel345}
\gamma=\lim_{n \rightarrow \infty} \left(1+\frac{1}{2}+\cdots +\frac{1}{n}-\log n\right)=0.57721566\ldots.
\end{equation}

In Section 2, we extend Euler's two zeta-function series for $\gamma$ to polylogarithm series for $\gamma(z)$ 
(Theorem~\ref{lilee}); we also show another way in which the function $\gamma(z)$ is related
to the extended polylogarithm (Theorem~\ref{veryi}). We give a variant of the series definition of $\gamma(z)$ which
is valid on a different domain (Theorem~\ref{fibbe233333}).  We prove two functional equations for $\gamma(z)$  
(Theorems~\ref{basiliki} and \ref{fibbe2}); the first is an inversion formula relating 
$\gamma(z)$ and $\gamma(1/z)$, and is proved 
using Jonqui\`{e}re's relation for the polylogarithm. 

In Section 3, we generalize one of Somos's constants, and show relations with the function $\gamma(z)$
(Theorem~\ref{sss111222}) and with an infinite nested radical due to Ramanujan (Corollary~\ref{atul_ram}). 
We also generalize Somos's 
quadratic recurrence sequence to cubic and other degrees
(Theorem~\ref{gngn}) and provide asymptotic estimates (Lemma~\ref{taylor9090} 
and Theorem~\ref{sondow11111}).

In Section 4, we calculate the value of $\gamma(z)$ at any root of unity (Theorem~\ref{ss1111}); the proof uses a result
based on  the Weierstrass product for the gamma function (Theorem~\ref{enana}).
Using the Kinkelin-Bendersky hyperfactorial $K$ function, in Section 5 we compute the derivative  $\gamma'(\omega)$ 
at a root of unity $\omega\neq 1$ (Theorem~\ref{0000ss1111}); the proof involves  a summation formula derived using 
the Barnes $G$ function (Theorem~\ref{e224nana}).
In particular, we show that $\gamma'(-1)$ 
involves the Glaisher-Kinkelin constant $A$ (Corollary~\ref{gammaprime1}), and that it is related to 
an infinite product (Example~\ref{eleven11}) essentially due to Borwein and Dykshoorn~\cite{bowr}.

Other infinite products occur in Corollary~\ref{three33},  Example~\ref{nine} (where we generalize an accelerated
product for pi~\cite{Sa111SS}), Section~3, and
 Example~\ref{severn7yes}. We evaluate some new double integrals in Examples~\ref{five} and \ref{seven}, 
Theorem~\ref{sumeee}, 
and Corollary~\ref{gammaprime1}.

\newpage

\noindent {\large {\bf 2. The Generalized-Euler-Constant Function $\gamma(z)$}}\label{twosection}

\vspace{.25in}

Denote by ${\mathbb Z}^+,$ $\R,$ and $\C$ the sets of positive integers, real numbers, and complex numbers,
respectively. If $z \in \C$ and $z\neq 0$,  define $$\log z=\ln |z|+i \,\Arg z \quad (-\pi < \Arg z \leq \pi).$$

The limit definition (\ref{defel345}) of Euler's constant is equivalent to the series formula
\begin{equation}\label{gammmmma}
\gamma=\sum_{n=1}^{\infty} \left(\frac{1}{n}-\log \frac{n+1}{n}\right)
\end{equation} 
(see \cite{sonqqq}).
The corresponding alternating series gives the ``alternating Euler constant" \cite{sonqqq}, \cite{quisond22} 
(see also \cite{quisond})
\begin{equation}\label{altegamma}
\log \frac{4}{\pi} = \sum_{n=1}^{\infty} (-1)^{n-1} \left(\frac{1}{n}-\log \frac{n+1}{n}\right)=0.24156447\ldots.
\end{equation}
The main subject of this paper is the following 
function $\gamma(z)$, which generalizes (\ref{gammmmma}) and (\ref{altegamma}).

\begin{definition}\label{defff1}
{\rm The {\it generalized-Euler-constant function} $\gamma(z)$ is defined when $|z|\leq 1$ by the power series
\begin{equation}\label{jamie1}
\gamma(z)=\sum_{n=1}^{\infty} z^{n-1} \left(\frac{1}{n}-\log\frac{n+1}{n}\right),
\end{equation}
which converges by comparison to series (\ref{gammmmma}) for $\gamma$.}
\end{definition}

\begin{example}\label{two}
{\rm In addition to $\gamma(1)=\gamma$ and $\gamma(-1)  =  \log\frac{4}{\pi},$
Definition~\ref{defff1} gives $\gamma(0)=1-\log2.$ At
$z=1/2$, the function takes the value
\begin{equation}\label{somnonss}
\gamma\left(\frac{1}{2}\right)=2\log\frac{2}{\sigma},
\end{equation}
where 
$$\sigma = \sqrt{1\sqrt{2\sqrt{3 \cdots}}}=1^{1/2} 2^{1/4} 3^{1/8} \cdots=\prod_{n=1}^{\infty} n^{1/2^n}=
1.66168794\ldots$$ is one of {\it Somos's quadratic recurrence constants}~\cite{somoss1}, 
\cite[Sequence A112302]{sloane222} 
(see also \cite{wwjam} and~\cite[p.\ 446]{FF}, where the notation $\gamma$ is used 
instead of $\sigma$). To see this, write (\ref{jamie1}) with $z=1/2$ as
\begin{eqnarray*}
\gamma\left(\frac{1}{2}\right) & = & 2 \sum_{n=1}^{\infty} \left( \frac{1}{2^{n}n}-
\frac{2\log (n+1)}{2^{n+1}}
+ \frac{\log n}{2^{n}}\right) \\
& = & 2\left(\log 2-2\log\sigma+\log\sigma\right) =2\log\frac{2}{\sigma}.
\end{eqnarray*}}
\end{example}

\vspace{.1in}

Euler gave two zeta-function series for his constant $\gamma$ (see, for example, \cite[equations 3.4(23) and 
3.4(151)]{sss}),
\begin{equation}\label{recover}
\gamma=\sum_{k=2}^{\infty} (-1)^k \frac{\zeta(k)}{k}=1-\sum_{k=2}^{\infty} \frac{\zeta(k)-1}{k}.
\end{equation}
We generalize them to polylogarithm series for the function $\gamma(z)$.

\begin{theorem}\label{lilee}
If $|z| \leq 1$ and if ${\rm Li}_k(z)$ denotes the {\rm polylogarithm} {\rm \cite[Section 1.11]{E}, \cite{quisond}, 
\cite{Sa111SS}}, defined for $k=2,3,\ldots$
by the convergent series 
$${\rm Li}_k(z)=\sum_{n=1}^{\infty} \frac{z^n}{n^k},$$
then 
\begin{equation}\label{kamalj}
z \gamma(z)=\sum_{k=2}^{\infty} (-1)^k\frac{{\rm Li}_k(z)}{k}.
\end{equation}
If in addition $z \neq 1$, then 
\begin{equation}\label{parebee}
z^2 \gamma(z)=z+(1-z)\log(1-z)-\sum_{k=2}^{\infty} \frac{{\rm Li}_k(z)-z}{k}.
\end{equation}
\end{theorem}
\noindent{\it Proof.} If $|z| \leq 1$, then Definition~\ref{defff1} and the expansion
\begin{equation}\label{formc}
\log(1-w)=-\sum_{m=1}^{\infty} \frac{w^m}{m} \quad (|w|\leq 1, w\neq 1)
\end{equation}
give
$$z\gamma(z)=z(1-\log 2)+\sum_{n=2}^{\infty} z^n \left[\frac{1}{n}-\log \left(1+\frac{1}{n}\right)\right]=z
\sum_{k=2}^{\infty} \frac{(-1)^k}{k} + \sum_{n=2}^{\infty} \sum_{k=2}^{\infty} (-1)^k \frac{z^n}{kn^k}.$$
It is easy to see that the double series converges absolutely, so we may 
reverse the order of summation, obtaining
$$z\gamma(z)=\sum_{k=2}^{\infty} \frac{(-1)^k}{k} \left( z+\sum_{n=2}^{\infty} \frac{z^n}{n^k}\right)=\sum_{k=2}^{\infty}
\frac{(-1)^k}{k} {\rm Li}_k(z).$$
This proves~(\ref{kamalj}).

Now take $z \neq 1$ with $|z|\leq 1.$ Definition~\ref{defff1} and formula (\ref{formc}) imply that 
$$z^2 \gamma(z)=z \sum_{n=1}^{\infty} \frac{z^n}{n}-\sum_{n=1}^{\infty} z^{n+1} \log \frac{n+1}{n}=-z
\log(1-z)-\sum_{n=2}^{\infty} z^n \log \frac{n}{n-1},$$
where we have re-indexed the last series. Thus, using (\ref{formc}) again, we can write
\begin{eqnarray*}
z^2 \gamma(z) & = & -z \log(1-z)+\log(1-z)+\sum_{n=1}^{\infty} \frac{z^n}{n}+\sum_{n=2}^{\infty} z^n
\log \frac{n-1}{n}\\
& = & (1-z)\log(1-z)+z+\sum_{n=2}^{\infty} z^n\left[\frac{1}{n}+\log\left(1-\frac{1}{n}\right)\right].
\end{eqnarray*}
The last series is equal to the absolutely convergent double series
$$-\sum_{n=2}^{\infty} \sum_{k=2}^{\infty} \frac{z^n}{kn^k}=-\sum_{k=2}^{\infty} \frac{1}{k} \sum_{n=2}^{\infty}
\frac{z^n}{n^k}=-\sum_{k=2}^{\infty} \frac{{\rm Li}_k(z)-z}{k},$$
and the proof of (\ref{parebee}) is complete.
\hfill $\Box$

\begin{example}\label{twoAAA}
{\rm Let $z=1$ in (\ref{kamalj}), and let $z$ tend to $1^{-}$ in (\ref{parebee}). 
Since $\gamma(1)=\gamma$ and 
${\rm Li}_k(1)=\zeta(k)$, we recover the two zeta series (\ref{recover}) for Euler's constant.

Now take $z=-1$, and substitute $\gamma(-1)=\log \frac{4}{\pi}$ and
${\rm Li}_k(-1)=(2^{1-k}-1)\zeta(k)$ (see \cite[Section 9.522]{soso00},
\cite{quisond}, \cite{Sa111SS}, \cite[Chapter 3]{spanoer}). Formula 
(\ref{kamalj}) and the first equality in (\ref{recover}) give the zeta series
(compare \cite[equation (4)]{sonqqq} and \cite[equation 3.4(25)]{sss})
$$\log \frac{4}{\pi}=\sum_{k=2}^{\infty} (-1)^k \frac{1-2^{1-k}}{k} \zeta(k)=\gamma-2\sum_{k=2}^{\infty} (-1)^k 
\frac{\zeta(k)}{2^k k}.$$
Simplification in formula (\ref{parebee}) yields a zeta series for the logarithm of pi:
$$\log \pi =1+\sum_{k=2}^{\infty} \frac{1-(1-2^{1-k})\zeta(k)}{k}.$$

Finally, set $z=1/2$. Using formula (\ref{somnonss}) for the Somos constant $\sigma$, and 
defining ${\rm Li}_1(z)=-\log(1-z)$, we can write formulas (\ref{kamalj}) and (\ref{parebee}) as
$$\log \sigma =\sum_{k=1}^{\infty} (-1)^{k-1}\frac{{\rm Li}_k(1/2)}{k}=\sum_{k=1}^{\infty} \frac{2{\rm Li}_k(1/2)-1}{k}.$$}
\end{example}

\begin{theorem}\label{journals}
The function $\gamma(z)$ is continuous on the closed unit disk $$D=\{z \in \C: |z|\leq 1\},$$
and holomorphic on the interior of $D$. 
However, the left-hand derivative at $z=1$ does not exist; more precisely,
\begin{equation}\label{heisnota}
\lim_{t \rightarrow 1^-}\frac{\gamma(1)-\gamma(t)}{1-t}=+\infty.
\end{equation}
In particular, the Taylor series expansion of $\gamma(z)$ at $z=0$ has radius of convergence $1$.
\end{theorem}
\noindent {\it Proof.} For $z \in D$,  series (\ref{jamie1}) is majorized by series (\ref{gammmmma}) 
for Euler's constant. Therefore, series (\ref{jamie1}) converges to $\gamma(z)$ uniformly on $D$. It follows that 
$\gamma(z)$ is continuous on $D$, and holomorphic on the interior of $D$. 

If $0 <t<1$, then 
\begin{eqnarray*}
{\gamma(1)-\gamma(t)} & = & \sum_{n=1}^{\infty} (1-t^{n-1})
\left[\frac{1}{n}-\log \left(1+\frac{1}{n}\right)\right]\\
& = & \sum_{n=1}^{\infty} (1-t^{n-1}) \left(\frac{1}{2n^2}-\frac{1}{3n^3}+\frac{1}{4n^4}- \cdots\right)
\geq \sum_{n=1}^N (1-t^{n-1}) \left(\frac{1}{2n^2}-\frac{1}{3n^3}\right)
\end{eqnarray*}
for $N=1,2,3,\ldots$. Multiplying by $(1-t)^{-1}$, and letting $t$ tend to $1^-$, we obtain 
$$\liminf_{t\rightarrow 1^-} \frac{\gamma(1)-\gamma(t)}{1-t} \geq \sum_{n=1}^N (n-1) \left(\frac{1}{2n^2}-
\frac{1}{3n^3}\right) \geq \sum_{n=1}^N \left(\frac{1}{2n}-\frac{5}{6n^2}\right).$$
The last sum tends to infinity with $N$, and the theorem follows.
\hfill $\Box$

\begin{figure}[t]
\setlength{\unitlength}{1cm} 
\begin{picture}(7,7)
\includegraphics{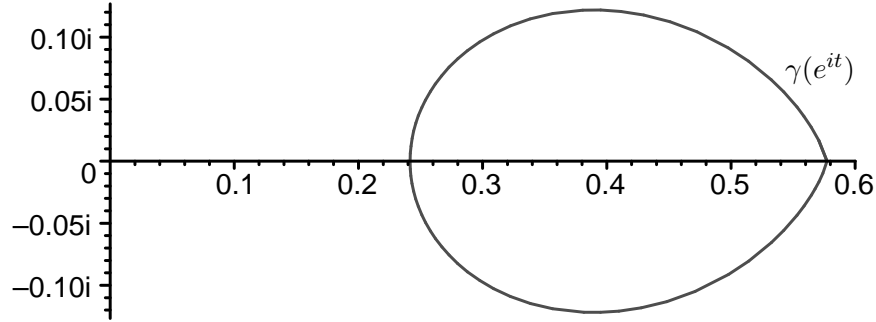} 
\end{picture}
\put (6,3.65){\small $\gamma(e^{it})$}
\caption{The image of the unit circle under the function $\gamma(z)$}\label{figg1}
\end{figure} \hfill

The image of the unit circle under the function $\gamma(z)$ is shown in Figure~\ref{figg1}.

In \cite{before, sonqqq, Sa111SS}, the first author represented Euler's constant $\gamma=\gamma(+1)$ and the
alternating Euler constant $\log \frac{4}{\pi}=\gamma(-1)$ by the double integrals
$$\gamma(\pm 1)=\int_0^1\int_0^1 \frac{1-x}{(1\mp xy)(-\log xy)}\,dx\,dy.$$
We extend this to integrals for the function $\gamma(z)$, and obtain its analytic continuation 
to the domain $\C-[1,\infty)$.

\begin{theorem}\label{identifu}
If $|z|\leq 1$, then 
\begin{equation}\label{sonora1}
\gamma(z)=\int_0^1\int_0^1 \frac{1-x}{(1-xyz)(-\log xy)}\,dx\,dy=\int_0^1 \frac{1-x+\log x}{(1-xz)\log x} \,dx.
\end{equation}
The integrals converge for all $z \in \C-(1,\infty)$, and provide the analytic continuation of 
the generalized-Euler-constant function $\gamma(z)$ for 
$z \in \C-[1,\infty).$ 
\end{theorem}
\noindent{\it Proof.}  By \cite[Theorem 4.1]{quisond},
if $z \in \C-(1,\infty)$ with $z \neq 0$, and if $\Re(s)>-2$ with $s\neq -1$, then for $u>0$
\begin{small}
\begin{equation}\label{exonoteso}
\int_0^1 \int_0^1 \frac{(1-x) (xy)^{u-1}}{(1-xyz) (-\log xy)^{-s}}\, dx \, dy  
 =  \Gamma(s+2) \left[
\Phi(z,s+2,u) + \frac{(1-z)\Phi(z,s+1,u)-u^{-s-1}}{z(s+1)}\right].  
\end{equation}
\end{small}\noindent  
Here $\Phi$ is the {\it Lerch transcendent} \cite[Section 1.11]{E},
\cite[Section 2]{quisond}, the analytic continuation of the series 
\begin{equation}\label{sstilef}
\Phi(z,s,u)=\sum_{n=0}^{\infty} \frac{z^n}{(n+u)^s},
\end{equation}
which converges for all complex $s$ and all $u>0$ when 
$|z|<1$. In particular, $\Phi(z,0,u)=1+z+z^2+\cdots=(1-z)^{-1}.$ It follows, letting $s$ tend to $-1$ 
in~(\ref{exonoteso}), that
\begin{equation}\label{iiitsiii}
\int_0^1 \int_0^1 \frac{(1-x) (xy)^{u-1}}{(1-xyz) (-\log xy)}\, dx \, dy  
 =  \Phi(z,1,u) + \frac{1-z}{z}\frac{\partial \Phi}{\partial s}(z,0,u)+\frac{\log u}{z}.
\end{equation}
Now set $u=1$ and let $I(z)$ denote the double integral
in (\ref{sonora1}). Using (\ref{sstilef}), and re-indexing, we obtain, when $0 < |z| <1$,
\begin{eqnarray}
I(z) & = & \sum_{n=1}^{\infty} \frac{z^{n-1}}{n} -(z^{-1}-1) \sum_{n=1}^{\infty} z^{n-1} \log n\nonumber \\
[-6pt] \label{ivtsi}\\
& = & \sum_{n=1}^{\infty} \frac{z^{n-1}}{n}-\sum_{n=1}^{\infty} z^{n-1} \log(n+1)+\sum_{n=1}^{\infty} z^{n-1}
\log n =\gamma(z),\nonumber
\end{eqnarray}
by Definition~\ref{defff1}.

This proves that $\gamma(z)=I(z)$ when $0 < |z|<1$. It follows, since the functions $\gamma(z)$ and $I(z)$ are 
continuous on $D-\{1\}$ (in fact, on $D$), that $\gamma(z)=I(z)$ also when $|z|=1\neq z$ and
when $z=0$. Finally, $\gamma(1)=\gamma=I(1)$, from \cite{sonqqq}.

It remains to show that $I(z)$ is equal to the single integral in (\ref{sonora1}). To see this, make the change of variables
$x=X/Y, y=Y$, and integrate with respect to $Y$.
\hfill $\Box$

\begin{figure}[t]
\setlength{\unitlength}{1cm} 
\begin{picture}(8,8)
\includegraphics{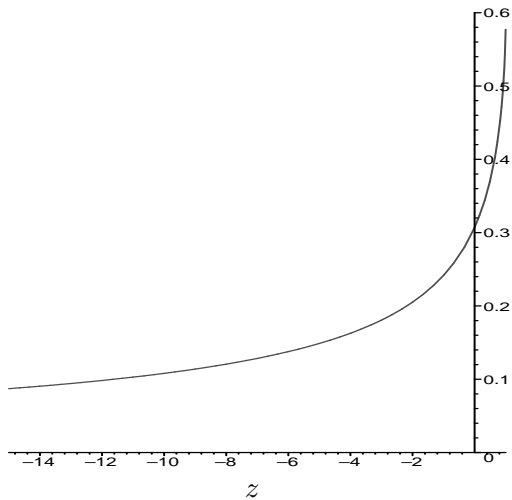} 
\end{picture}
\put (0.3,-0.1){\small $z$}
\caption{The function $\gamma(z)$ for real $z$}\label{figg2}
\end{figure} \hfill

The graph of the function $\gamma(z)$ for $z \in (-\infty,1]$ is shown in Figure~\ref{figg2}. 
Its properties are easily verified from (\ref{sonora1}) and (\ref{heisnota}). 
Namely, {\it for $z$ real, the graph of $\gamma(z)$ is positive, increasing, and concave upward; it is 
asymptotic to the negative real axis, that is}, $\lim_{z \rightarrow \infty} \gamma(-z)=0$; {\it and the tangent line
at the point $(1,\gamma)$ is vertical.}

\begin{example}\label{three} 
{\rm Using Theorem~\ref{identifu} and formula (\ref{somnonss}) for Somos's constant $\sigma$, we recover 
the evaluations from \cite{quisond}
$$\int_0^1 \int_0^1 \frac{x}{(2-xy)(-\log xy)}\,dx\,dy=\int_0^1 \frac{1-x}{(2-x)(-\log x)} \, dx=\log\sigma.$$
To see this, set $z=1/2$ in (\ref{sonora1}), and replace $\gamma(1/2)$ with $2\log (2/\sigma)$. Now multiply 
 by $1/2$, and subtract the result from the equations
$$\int_0^1 \int_0^1 \frac{1}{(2-xy)(-\log xy)}\,dx\,dy=\int_0^1 \frac{1}{2-x}\,dx=\log 2,$$
which are easily verified.}
\end{example}

\vspace{.1in}

The function $\gamma(z)$ is also related to the polylogarithm in a different way from that in Theorem~\ref{lilee}.

\begin{theorem}\label{veryi}
For all $z \in \C-[1,\infty)$ the relation
\begin{equation}\label{1par}
z^2\gamma(z)=(1-z){\rm Li}_0'(z)-z\log(1-z)
\end{equation}
holds. Here the prime $'$ denotes $\partial/\partial s$, and ${\rm Li}_s(z)$ is the 
{\rm extended polylogarithm}, the analytic continuation {\rm \cite[Section 1.11]{E}}, 
{\rm \cite[Section~2]{quisond},}
{\rm \cite[Section 5]{Sa111SS}} of the series 
\begin{equation}\label{posesfores111}
{\rm Li}_s(z)=\sum_{n=1}^{\infty} \frac{z^n}{n^s}.
\end{equation}
If $|z|<1$, then the series converges for any $s \in \C$, and the relation can be written
\begin{equation}\label{viitsi}
z\gamma(z)=-\log(1-z)-(1-z)\sum_{n=1}^{\infty} z^{n-1} \log n.
\end{equation}
\end{theorem}
We give two proofs.

\vspace{.1in}

\noindent{\it Proof 1.} Multiply the last two equations in (\ref{ivtsi}) by $z$, with $|z|<1$. Using (\ref{formc}), formula 
(\ref{viitsi}) follows. Now multiply (\ref{viitsi}) by $z$. Using (\ref{posesfores111}), we obtain relation (\ref{1par}), and
the theorem follows by analytic continuation. \hfill $\Box$

\vspace{.1in}

\noindent{\it Proof 2.} Setting $u=1$ in (\ref{iiitsiii}), we may replace the integral with $\gamma(z)$, by 
Theorem~\ref{identifu}. Now multiply the equation by $z^2$. Using the formulas $-\log(1-z)=z\Phi(z,1,1)$ (from
(\ref{formc}) and (\ref{sstilef})) and 
\begin{equation}\label{ss11nn11}
{\rm Li}_s(z)=z\Phi(z,s,1)
\end{equation}
(from (\ref{posesfores111}) and (\ref{sstilef})), relation (\ref{1par}) follows. 
If $|z|<1$, substituting (\ref{posesfores111}) in (\ref{1par}) yields (\ref{viitsi}).
\hfill $\Box$

\vspace{.1in}

An application of Theorem~\ref{veryi} is a formula for $\gamma(z)$ involving a series which converges if
$\Re(z) <1/2$, when $\left| \frac{-z}{1-z}\right| <1$.

\begin{theorem}\label{fibbe233333}
If $\Re(z) < 1/2$, then
\begin{equation}\label{jamie}
z\gamma(z)=-\log(1-z)+\sum_{n=1}^{\infty} \left(\frac{-z}{1-z}\right)^{n} \sum_{k=0}^n (-1)^{k+1} 
{n \choose k} \log (k+1).
\end{equation}
Equivalently, if $|w| <1$, then
\begin{equation}\label{jamieee}
w\gamma\left(\frac{-w}{1-w}\right)=-(1-w)\log(1-w)+(1-w)\sum_{n=1}^{\infty} w^n \sum_{k=0}^n 
(-1)^{k} {n \choose k} \log(k+1).
\end{equation}
\end{theorem}
\noindent {\it Proof.}
By a special case of  \cite[Corollary 5.1]{quisond}, if $\Re(z) <1/2$, then
the double sum in (\ref{jamie}) is equal to the product $(1-z)\frac{\partial \Phi}{\partial s}(z,0,1)$.
Using (\ref{ss11nn11}) and (\ref{1par}), the first statement follows. To prove the second, note that if $|w|<1$ and
$z=\frac{-w}{1-w}$, then $\Re(z)<1/2$, and (\ref{jamie}) implies (\ref{jamieee}). \hfill $\Box$

\begin{corollary}\label{three33}
For $m=1,2,3,\ldots,$ define the infinite product
$$P_m=\prod_{n=1}^{\infty} \left(\prod_{k=0}^n (k+1)^{(-1)^{k+1}{n \choose k}}\right)^{\left(\frac{m}{m+1}\right)^n}.$$
Then the product converges, and
$$\gamma(-m)=\frac{1}{m} \log \frac{m+1}{P_m}.$$
\end{corollary}
\noindent {\it Proof.} When $z=-m$, the series in (\ref{jamie}) is equal
to $\log P_m$. \hfill $\Box$

\vspace{.1in}

\begin{example}\label{nine}
{\rm Let $m=1$. The product $P_1$ is the acceleration of Wallis's product for 
$\pi/2$ in~\cite{Sa111SS} (see also \cite{quisond}):
$$P_1=\left(\frac{2}{1}\right)^{1/2} \left(\frac{2^2}{1\cdot 3}\right)^{1/4} \left(\frac{2^3\cdot 4}{1 \cdot 3^3}\right)^{1/8}
\left( \frac{2^4 \cdot 4^4}{1 \cdot 3^6 \cdot 5}\right)^{1/16} \cdots =\frac{\pi}{2}.$$
Thus
$$\gamma(-1)=\log\frac{2}{P_1}=\log \frac{4}{\pi},$$
confirming the value of $\gamma(-1)$ in Example \ref{two}.

With $m=2$, we get 
$$\gamma(-2)=\frac{1}{2} \log \frac{3}{P_2},$$
where 
$$P_2=\left(\frac{2}{1}\right)^{2/3} \left(\frac{2^2}{1\cdot 3}\right)^{4/9} \left(\frac{2^3\cdot 4}{1 \cdot 3^3}\right)^{8/27}
\left( \frac{2^4 \cdot 4^4}{1 \cdot 3^6 \cdot 5}\right)^{16/81} \cdots .$$}
\end{example}

\vspace{.1in}

We now prove a functional equation for the generalized-Euler-constant function 
 which expresses $\gamma(1/z)$ in terms of $\gamma(z)$.

\begin{theorem} \label{basiliki}
The following inversion formula holds for all $z \in \C-[0,\infty)$ with $\Im(z) \geq 0$:
$$\gamma\left(\frac{1}{z}\right)=z^3\gamma(z)-\pi i z+z^2\log(1-z)-z \log \frac{z-1}{z}+z(1-z)\left[\gamma+\log 2\pi
+\frac{\pi i}{2}+\psi \left( \frac{\log z}{2\pi i}\right)\right],$$
where $\psi(z)$ is the {\rm digamma function} $$\psi(z)=\frac{d[\log\Gamma(z)]}{dz}=\frac{\Gamma'(z)}{\Gamma(z)}.$$
\end{theorem}
\noindent {\it Proof.}
The extended polylogarithm ${\rm Li}_s(z)$ (see Theorem~\ref{veryi}) satisfies {\it Jonqui\`{e}re's relation}
\cite[Section 1.11]{E}, \cite{jonqqq}
\begin{equation}\label{a24}
{\rm Li}_s(z)+e^{\pi i s}{\rm Li}_s(1/z)=\frac{(2\pi)^s e^{\pi i s/2}}{\Gamma(s)} \,\zeta\!\left(1-s, \frac{\log z}{2\pi i}\right),
\end{equation}
where $\zeta(s,w)$ is the {\it Hurwitz} (or {\it generalized}) {\it zeta function}~\cite[Section 1.10]{E}, 
the analytic continuation of the series
$$\zeta(s,w)=\sum_{n=0}^{\infty} \frac{1}{(n+w)^s},$$
which converges if $\Re(s)>1$ and $w \in \C-\{0,-1,-2,\ldots\}.$ If $|z|<1$ and $s=0$, the sum of the series 
(\ref{posesfores111}) for ${\rm Li}_s(z)$ is ${\rm Li}_0(z)=z(1-z)^{-1},$
which by analytic continuation holds for all complex $z \neq 1.$ It follows that at $s=0$ the left side of (\ref{a24}) is 
equal to $-1$. Therefore, the derivative of (\ref{a24}) with respect to $s$ at $s=0$ is
$${\rm Li}_0'(z)+{\rm Li}_0'(1/z)
+\frac{\pi i}{z-1}=\lim_{s\rightarrow 0} \frac{1}{s}
\left[\frac{(2\pi)^s e^{\pi i s/2}}{\Gamma(s)}\,\zeta\! \left(1-s,\frac{\log z}{2\pi i}\right)-(-1)\right].$$
To compute the limit, we use the Taylor series for 
$$(2\pi)^s e^{\pi i s/2}=\exp\left[\left(\log 2\pi +\frac{\pi i}{2}\right)s\right]$$
together with the estimates \cite[Equation 43:6:1]{spanoer}
$$\frac{1}{\Gamma(s)} = s+\gamma s^2 +O(s^3)$$
and \cite[Section 1.11]{E}
$$\zeta(1-s,w)=-\frac{1}{s}-\psi(w)+O(s) \quad (\Re(w) >0),$$
which are valid for $s$ tending to $0$. If $\Im(z) >0$, the result is
$${\rm Li}_0'(z)+{\rm Li}_0'(1/z)
+\frac{\pi i}{z-1}=-\gamma-\log 2\pi -\frac{\pi i}{2}
-\psi\left(\frac{\log z}{2\pi i}\right).$$
For $z \in \C-[0,\infty)$, we may use (\ref{1par}) to replace ${\rm Li}_0'(z)$ and ${\rm Li}_0'(1/z)$ 
with expressions for them involving $\gamma(z)$ and $\gamma(1/z)$, respectively. 
Solving for $\gamma(1/z)$, we arrive at the inversion formula. This proves it when $\Im(z)>0$.
Since each term in the formula is continuous on the set $\{z\in \C| \, \Re(z)<0, \, \Im(z)\geq 0\}$, the theorem follows.
\hfill $\Box$


\begin{example}\label{ten}
{\rm Take $z=-1$. Using the values $\log (-1)=\pi i$ and $\psi(1/2)=-\gamma-\log 4$ (from Gauss's formula
\cite[equation 1.2(47)]{sss}
\begin{equation}\label{gauss}
\psi\left(\frac{j}{q}\right)  =  -\gamma-\frac{\pi}{2}\cot \frac{j\pi}{q}-\log q + \sum_{k=1}^{q-1} \cos
\frac{2kj\pi}{q}\log\left(2 \sin \frac{k\pi}{q}\right),
\end{equation}
where $0 < j < q$), we again obtain $\gamma(-1)=\log\frac{4}{\pi}.$ For a related application of Theorem~\ref{basiliki}, 
see the proof of Corollary~\ref{gammaprime1}.}
\end{example}


\begin{remark}
{\rm Setting $z=1/w$ in Theorem~\ref{basiliki}, we obtain an inversion formula valid for $w \in \C-[0,\infty)$ with 
$\Im (w) \leq 0$.}
\end{remark}

The next theorem gives a second functional equation for the function $\gamma(z)$. The equation relates the 
quantities $\gamma(z), \gamma(-z)$, and $\gamma(z^2)$.

\begin{theorem}\label{fibbe2}
The following reflection formula holds
for all $z\in \C-((-\infty,-1]\cup [1,\infty))$:
$$z(1+z)\gamma(z)+z(1-z)\gamma(-z)=2z^3 \gamma(z^2)-2z\log 2+(1+z)\log (1+z)-(1-z)\log(1-z).$$
\end{theorem}
\noindent {\it Proof.}
If $|z|<1$, then using (\ref{posesfores111}) we see that
$$\frac{1}{2}\left({\rm Li}_0'(z)+{\rm Li}_0'(-z)\right)=-\sum_{n=1}^{\infty} z^{2n} \log 2n =
\frac{z^2\log 2}{z^2-1}+{\rm Li}_0'(z^2),$$
where the prime $'$ denotes $\partial/\partial s$. The relation between $\gamma(z)$ and ${\rm Li}'_0(z)$ in
Theorem~\ref{veryi} then yields the desired formula, and the result follows by analytic continuation.
\hfill $\Box$

\vspace{.1in}

\begin{remark} {\rm Theorem~\ref{fibbe2} can be generalized, as follows. Given an integer $q>1$, we choose
a $q$th root of unity $\omega \neq 1$, and obtain the average
$$\frac{1}{q} \sum_{j=0}^{q-1} {\rm Li}_0'(\omega^jz)=\frac{z^q \log q}{z^q-1}+{\rm Li}_0'(z^q).$$
Theorem~\ref{veryi} then translates this into a formula relating $\gamma(z), \gamma(\omega z), \ldots,
\gamma(\omega^{q-1}z)$ and $\gamma(z^q)$, valid for all $z \in \C$ such that $\omega^j z \not \in [1,\infty)$ for 
$j=0,1,\ldots,q-1$.}
\end{remark}

\vspace{.25in}

\noindent {\large {\bf 3. A Generalization of Somos's Quadratic Recurrence Constant}}\label{threesection}

\vspace{.25in}

We begin this section by  generalizing both Somos's quadratic recurrence constant $\sigma$
and its relation (\ref{somnonss}) with the function $\gamma(z)$.
By convention, if $a\geq 0$ and $t > 0$,  we assume that $\sqrt[t]{a}=a^{1/t} \geq 0.$

\begin{definition}\label{controldd}
{\rm For $t >1$, the {\it generalized Somos constant} $\sigma_t$ is given by
$$\sigma_t=\sqrt[t]{1\sqrt[t]{2\sqrt[t]{3\cdots}}}=1^{1/t} 2^{1/t^2}3^{1/t^3} \cdots=\prod_{n=1}^{\infty} n^{1/t^n}.$$}
\end{definition}

The convergence of the infinite product for $\sigma_t$  follows from the convergence of the series
\begin{equation}\label{angolabro}
\log \sigma_t=\sum_{n=1}^{\infty} \frac{\log n}{t^n}=-{\rm Li}_0'\left(\frac{1}{t}\right)=
-\frac{1}{t} \frac{\partial \Phi}{\partial s}\left(\frac{1}{t},0,1\right),
\end{equation}
where $\Phi$ is the  Lerch transcendent {\rm (\ref{sstilef})}.
Note that for $t=2$  we get Somos's constant $\sigma=\sigma_2.$

The following result relates the generalized Somos constant $\sigma_t$ to the function $\gamma(z)$, essentially
generalizing Example~\ref{two}.

\begin{theorem} \label{sss111222}
For $t >1$, the generalized Euler constant $\gamma(1/t)$ and the generalized Somos
constant $\sigma_t$ satisfy the relation
\begin{equation}\label{ss1100a}
\gamma\left(\frac{1}{t}\right)=t \log \frac{t}{(t-1) \sigma_t^{t-1}}.
\end{equation}
In particular, 
\begin{equation}\label{somosuytr}
\lim_{t \rightarrow 0^+} t \sigma_{t+1}^t=e^{-\gamma}.
\end{equation}
\end{theorem}
\noindent {\it Proof.} To prove (\ref{ss1100a}), set $z=1/t$ in (\ref{viitsi}) and use (\ref{angolabro}). 
(The case $t=2$ is proved in Example~\ref{two}.)
To prove (\ref{somosuytr}), replace $t$ with $t+1$ in (\ref{ss1100a}), let $t$ tend to $0^+$, and use 
$\lim_{x \rightarrow 1^-} \gamma(x)=\gamma(1)=\gamma$.~\hfill $\Box$

\begin{remark}\label{tworrr}
{\rm From Theorem~\ref{sss111222} and the reflection formula in Theorem~\ref{fibbe2},
one can also express $\gamma(-1/t)$ 
in terms of $\sigma_t$ and $\sigma_{t^2}$.}
\end{remark}

The next theorem generalizes a result due to Somos~\cite{somoss1} (see \cite[p.\ 446]{FF} and \cite{wwjam}).

\begin{theorem}\label{gngn}
Fix $t>1$, and define the sequence $(g_{n,t})_{n \geq 0}$ via the generalized Somos recurrence
\begin{equation}\label{somoslae}
g_{0,t}=1,  \;\; g_{n,t}=ng_{n-1,t}^t \;\;\; (n\geq 1).
\end{equation}
Then we have the explicit solution
$$g_{n,t}=\sigma_t^{t^n} \exp\left[\frac{1}{t} \frac{\partial \Phi}{\partial s}\left(\frac{1}{t},0,n+1\right)\right]
=\sigma_t^{t^n}  \prod_{m=1}^{\infty} (m+n)^{-1/t^m}.$$
\end{theorem}
\noindent {\it Proof.} This follows by induction on $n.$ \hfill $\Box$

\vspace{.1in}

For $t=2$, equation (\ref{somoslae}) is Somos's quadratic recurrence $g_n=n g_{n-1}^2$. 
Somos \cite{somoss1}, \cite[Sequence A116603]{sloane222} (see also \cite{FF} and \cite{wwjam}) 
gave the following asymptotic formula for $g_n=g_{n,2}$ as $n$ tends to infinity:
\begin{equation}\label{tookmybooks}
g_n \sim \sigma^{2^n}(n+2-n^{-1}+4n^{-2}-21n^{-3}+138n^{-4}-1091n^{-5}+\ldots)^{-1}.
\end{equation}
We extend this to an asymptotic formula for $g_{n,t}$ given any fixed $t>1$. 

\begin{lemma}\label{taylor9090}
For $t>1$, let $(g_{n,t})_{n\geq 0}$ be the sequence defined by the recurrence {\rm (\ref{somoslae})}.
For $x\in[0,\infty)$, define $f_t(x) \in [1,\infty)$ by the infinite product
$$f_t(x)=\prod_{m=1}^{\infty}(1+mx)^{1/t^m}.$$
Then for each $N\geq 2$ and $n\geq 1$, there exists a positive number $\mu=\mu(N,n,t)$ such that 
\begin{equation}\label{accident}
g_{n,t}=\sigma_t^{t^n} n^{-1/(t-1)} \left[1+\sum_{k=1}^{N-1} \frac{1}{n^k k!}
 \frac{\partial^k f_t}{\partial x^k}(0)+\frac{1}{n^N N!}
 \frac{\partial^N f_t}{\partial x^N}(\mu)\right]^{-1}.
\end{equation}
Moreover, for fixed $N\geq 2$ and fixed $t>1$, 
$$\lim_{n\rightarrow \infty} \mu(N,n,t)=0, \quad 
\lim_{n\rightarrow \infty} \frac{\partial^N f_t}{\partial x^N}(\mu(N,n,t))=\frac{\partial^N f_t}{\partial x^N}(0).$$
\end{lemma}
\noindent {\it Proof.} Since $t>1$ and $x\geq 0$, the product for $f_t(x)$ converges.
It follows from Theorem~\ref{gngn} and the identity $\sum_{m=1}^{\infty}t^{-m}=1/(t-1)$ that
$$g_{n,t}=\sigma_t^{t^n} n^{-1/(t-1)} f_t\left(\frac{1}{n}\right)^{-1}.$$
Thus, by Taylor's theorem with remainder, it suffices to show that $f_t(x)$ is infinitely differentiable with respect
to $x$ for $x \in [0,\infty)$, the derivatives at $x=0$ being right-sided. That in turn follows (see, for example, 
\cite[p.\ 342, Theorem~4]{knopp}) by noting that for $k \geq 1$ the $k$th termwise derivative of the series for $\log f_t(x)$,
namely, $$(-1)^{k-1}(k-1)!
\sum_{m=1}^{\infty} \frac{m^k}{t^m(1+mx)^k},$$ is uniformly convergent on $[0,\infty)$, as it is majorized by
$(k-1)! \sum_{m=1}^{\infty} m^k t^{-m}$. \hfill $\Box$

\begin{theorem}\label{sondow11111}
For fixed $t>1$, the sequence $(g_{n,t})_{n\geq 0}$ satisfies the following asymptotic condition
as $n$ tends to infinity:
\begin{small}
$$g_{n,t}=\sigma_t^{t^n} n^{-1/(t-1)} \left[1+\frac{t}{(1-t)^2n}-\frac{t(t^2-t-1)}{2(1-t)^4n^2}+
\frac{t(2t^4+t^3-11t^2+7t+2)}{6(1-t)^6n^3} +O\left(\frac{1}{n^4}\right)
\right]^{-1}.$$
\end{small}
\end{theorem}
\noindent {\it Proof.} By Lemma~\ref{taylor9090}, we only have to compute $\frac{\partial^k f_t}{\partial x^k}(0)$.
For $k=1,2,\ldots$ define $\phi_k(t)$ by 
$$\phi_k(t)=\frac{\partial^k \log f_t}{\partial x^k}(0)=(-1)^{k-1}(k-1)!\sum_{m=1}^{\infty} \frac{m^k}{t^m} 
=(-1)^{k-1}(k-1)! {\rm Li}_{-k}\left(\frac{1}{t}\right).$$
Using the identity
$$\phi_k(t)=(k-1)t\,\frac{\partial \phi_{k-1}}{\partial t}(t) \quad (k>1),$$
we compute that 
$$\phi_1(t)  =  \frac{t}{(1-t)^2}, \quad \phi_2(t)  =  \frac{t(1+t)}{(1-t)^3},\quad 
\phi_3(t)  =  \frac{2t(t^2+4t+1)}{(1-t)^4}.$$

Since $f_t=e^{\log f_t}$ and $f_t(0)=1$, Fa\'{a} di Bruno's formula for the $k$th derivative 
of the composition of two functions
(see, for example, \cite[p.\ 21]{soso00}) yields
$$\frac{\partial^k f_t}{\partial x^k}(0)= \sum \frac{k!}{m_1!\cdots m_k!}
\prod_{j=1}^k \left(\frac{\phi_j(t)}{j!}\right)^{m_j},$$
where the sum is over all $k$-tuples $(m_1,\ldots,m_k)$ of non-negative integers such that 
$$m_1+2m_2+3m_3+\ldots+km_k=k.$$ In particular, 
$$\frac{\partial f_t}{\partial x}(0)  =  \phi_1(t), \quad
\frac{\partial^2 f_t}{\partial x^2}(0)  =  \phi_1(t)^2+\phi_2(t), \quad
\frac{\partial^3 f_t}{\partial x^3}(0)  =  \phi_1(t)^3+3\phi_1(t)\phi_2(t)+\phi_3(t).$$
Applying  (\ref{accident}) with $N=4$, the theorem follows.
\hfill $\Box$

\vspace{.1in}

\begin{example}{\rm
Taking $t=2$ gives the first four terms of the asymptotic formula (\ref{tookmybooks}) 
for Somos's quadratic recurrence sequence 
$g_n=g_{n,2}=ng_{n-1,2}^2=1,2,12,576,1658880,\ldots$ \cite[Sequence A052129]{sloane222}. 
With $t=3$ we get an asymptotic formula \cite[Sequences A123853 and A123854]{sloane222} 
for the cubic recurrence sequence $g_{n,3}=n g_{n-1,3}^3=1,2,24,55296,845378412871680,\ldots$ 
\cite[Sequence A123851]{sloane222}, namely,
$$g_{n,3} = \sigma_3^{3^n} n^{-1/2}\left[1+\frac{3}{4n}-\frac{15}{32n^2}+\frac{113}{128n^3}+O\left(\frac{1}{n^4}\right)
\right]^{-1}\quad (n \rightarrow \infty),$$
where $\sigma_3$ is the cubic recurrence constant \cite[Sequence A123852]{sloane222}
$$\sigma_3 = \sqrt[3]{1\sqrt[3]{2\sqrt[3]{3 \cdots}}}=1^{1/3} 2^{1/9} 3^{1/27} \cdots = 1.15636268\ldots.$$}
\end{example}

\vspace{.1in}

The generalized Somos constants $\sigma_t$ 
are connected with some of Ramanujan's infinite nested radicals and with 
the Vijayaraghavan-Herschfeld convergence criterion for them, as follows. Fix $t >1$.
Given $a_n =a_n(t) \geq 0$ for $n=1,2,\ldots,$ 
define $b_n=b_n(t)$ by
\begin{equation}\label{baktiar}
b_n = \sqrt[t]{a_1+\sqrt[t]{a_2+\sqrt[t]{a_3+\cdots+\sqrt[t]{a_{n-1}+\sqrt[t]{a_n}}}}}.
\end{equation}
If the sequence 
$(b_n)_{n\geq 1}$ converges, we call the limit $b_{\infty}=b_{\infty}(t)$, and we define
\begin{equation}\label{nested11}
\sqrt[t]{a_1+\sqrt[t]{a_2+\sqrt[t]{a_3+\cdots}}}=\lim_{n \rightarrow \infty} b_n=b_{\infty}.
\end{equation}

The following theorem is a special case of a result stated by Herschfeld~\cite[Theorem III]{HERR}. For $t=2$, it was first
proved by Vijayaraghavan (in Appendix I of \cite[p.\ 348]{dajur1}), 
and independently by Herschfeld.
(The proof for $t>1$ is similar to that for $t=2$, and thus is omitted by both authors.)

\begin{theorem}\label{twelve12} {\bf (Vijayaraghavan-Herschfeld)}
Fix $t >1$, and suppose that $a_n=a_n(t)\geq 0$ for $n\geq 1$. 
The sequence $(b_n)_{n \geq 1}$ defined by {\rm (\ref{baktiar})} converges if and only if 
$$\limsup_{n \rightarrow \infty} a_n^{1/t^n} < \infty.$$
\end{theorem}
Note that this inequality is equivalent to  the sequence $(a_n^{1/t^n})_{n \geq 1}$ 
being bounded. 

We now give the promised connection between the generalized Somos constant $\sigma_t$ and Ramanujan's infinite
radical. The latter is defined by 
\begin{equation}\label{nested22}
\sqrt[t]{1+2\sqrt[t]{1+3\sqrt[t]{1+\cdots}}}=\lim_{n \rightarrow \infty} B_n(t)=B_{\infty}(t),
\end{equation}
where 
$$B_n(t)=\sqrt[t]{1+2\sqrt[t]{1+3\sqrt[t]{1+\cdots+(n-1)\sqrt[t]{1+n\sqrt[t]{1}}}}}.$$

\begin{corollary}\label{atul_ram}
Fix $t>1$.
If  $$a_n=a_n(t)=1^{t^n} \cdot 2^{t^{n-1}} \cdot 3^{t^{n-2}} \cdots n^{t^1}$$
for $n=1,2,\ldots$, then 
\begin{equation}\label{somos12300}
\lim_{n \rightarrow \infty} a_n^{1/t^n}=\sigma_t^t < \infty,
\end{equation}
and the infinite nested radicals $b_{\infty}(t)$ and $B_{\infty}(t)$ in {\rm (\ref{nested11})} and {\rm 
(\ref{nested22})} converge and are equal.
\end{corollary}
\noindent {\it Proof.} The equality~(\ref{somos12300}) follows immediately from Definition~\ref{controldd}.
Since $\sigma_t^t < \infty$, Theorem~\ref{twelve12} implies that
the infinite nested radical 
$$\sqrt[t]{a_1+\sqrt[t]{a_2+\sqrt[t]{a_3+\cdots}}}=
\sqrt[t]{1^{t^1}+\sqrt[t]{1^{t^2} \cdot 2^{t^1}+\sqrt[t]{1^{t^3}\cdot 2^{t^2} \cdot 3^{t^1}+\cdots}}}=\lim_{n \rightarrow \infty}
b_n(t)=b_{\infty}(t)$$
converges. Now notice that $b_n(t)=B_n(t)$ 
for $n=1,2,\ldots$. Therefore, the infinite nested radical $B_{\infty}(t)$ 
also converges, and $b_{\infty}(t)=B_{\infty}(t).$ This proves the corollary.
\hfill $\Box$

\begin{remark} {\rm The only known value of $B_{\infty}(t)$ 
is $B_{\infty}(2)=3,$  that is,
$$\sqrt{1+2\sqrt{1+3\sqrt{1+\cdots}}}=3.$$
This formula was discovered by 
Ramanujan \cite[Question 289 and Solution, p.\ 323]{dajur1}, but his proof
is incomplete. Vijayaraghavan \cite[p. 348]{dajur1} and Herschfeld~\cite[pp.\ 420-421]{HERR} each completed 
Ramanujan's 
proof. The result also appeared as Putnam problem A6 in 1966 (see \cite[pp.\ 5 and 52]{alexxx111}).}
\end{remark}

\vspace{.25in}

\noindent {\large {\bf 4. Calculation of $\gamma(z)$ at Roots of Unity}}\label{onesection}

\vspace{.25in}

The purpose of this section is to calculate the value of $\gamma(z)$ at roots of unity. We use
the following  summation formula, which has some resemblance to Problem 24 in~\cite[p.\ 71]{sss}.

\begin{theorem}\label{enana}
Let $(\theta_n)_{n \geq 1}$ be a periodic sequence in $\C$, and let $q$ be a positive multiple of its period. 
If $\sum_{j=1}^q \theta_j=0$, then 
$$\sum_{n=1}^{\infty}   \theta_n \log\frac{n+1}{n}=
\sum_{j=1}^q \theta_j \log
\frac{\Gamma\left(\frac{j}{q}\right)}{\Gamma\left(\frac{j+1}{q}\right)}.$$
\end{theorem}
\noindent {\it Proof.}
The Weierstrass product for the gamma function is \cite[Section 8.322]{soso00}
$$\Gamma(z)=\frac{e^{-\gamma z}}{z} \prod_{n=1}^{\infty} \frac{e^{z/n}}{1+\frac{z}{n}}$$
for $z \in \C-\{0,-1,-2,\ldots\}.$  It follows that
for $j=1,\ldots,q$,
\begin{equation}\label{callclm}
\log \Gamma \left(\frac{j}{q}\right)-\log\Gamma\left(\frac{j+1}{q}\right)=\frac{\gamma}{q}
+\log \frac{j+1}{j}+\sum_{n=1}^{\infty} \left(-\frac{1}{qn}+\log \frac{qn+j+1}{qn+j}\right).
\end{equation}
Now multiply by $\theta_j$, and sum from $j=1$ to $q$.  Using $\sum_{j=1}^q \theta_j=0$ and
$\theta_j=\theta_{qn+j}$, we obtain
$$\sum_{j=1}^q \theta_j \log \frac{\Gamma\left(\frac{j}{q}\right)}{\Gamma\left(\frac{j+1}{q}\right)}=
\sum_{j=1}^q \theta_j \log \frac{j+1}{j}+
\sum_{n=1}^{\infty} \sum_{j=1}^q \theta_j \log \frac{qn+j+1}{qn+j}=\sum_{n=1}^{\infty} \theta_n \log \frac{n+1}{n},$$
where the last series converges by Dirichlet's test. This proves the theorem.
\hfill $\Box$

\vspace{.1in} 

\begin{example} \label{severn7yes}
{\rm Take $\theta_n=(-1)^{n-1}$ and $q=2$, and  exponentiate the series. Using the identities
(valid for $z \in \C-\{0,-1,-2,\ldots\}$) 
\begin{equation}\label{hajaja}
\Gamma(z+1)=z\Gamma(z)
\end{equation}
and
\begin{equation}\label{77sss}
\Gamma(z)\Gamma(1-z)=\frac{\pi}{\sin \pi z}
\end{equation}
to compute $\Gamma(1/2)\Gamma(3/2)=\pi/2$, we recover Wallis's product for pi,
which can be written
$$\prod_{n=1}^{\infty} \frac{4n^2}{4n^2-1}=\frac{\pi}{2}.$$

Now take $\theta_n=i^{n-1}$ and $q=4$. Exponentiating the real and imaginary parts of the series, and putting $z=1/4$
in (\ref{77sss}), we get the pair of products
$$ \prod_{n=1}^{\infty} \left(\frac{2n}{2n-1}\right)^{(-1)^{n-1}}=\frac{\Gamma(1/4)^2}{\pi^{3/2}\sqrt{2}},
\quad \prod_{n=1}^{\infty} \left(\frac{2n+1}{2n}\right)^{(-1)^{n-1}}=\frac{\Gamma(1/4)^2}{4\sqrt{2\pi}}.$$
Multiplying them together, or dividing the second by the first, we obtain the formulas
$$\prod_{n=1}^{\infty} \left(\frac{2n+1}{2n-1}\right)^{(-1)^{n-1}}=\frac{\Gamma(1/4)^4}{8\pi^2}, \quad
\prod_{n=1}^{\infty} \left(\frac{4n^2-1}{4n^2}\right)^{(-1)^{n-1}}=\frac{\pi}{4}.$$
The last may be a new product for pi. If we multiply it by Wallis's product, the factors with odd
$n$ cancel, and taking the square root gives
$$\prod_{n=1}^{\infty} \frac{16n^2}{16n^2-1} =\frac{\pi}{2\sqrt{2}}.$$
These products are perhaps known.
}
\end{example}

The next result gives two formulas for the value of $\gamma(z)$ at any root of unity $\omega \neq 1$.

\begin{theorem}\label{ss1111}
Let $p$ and $q$ be relatively prime positive integers. If $\omega=e^{i\pi p/q} \neq 1$, then
$$\gamma(\omega)+\frac{\log(1-\omega)}{\omega}=\sum_{j=1}^q \omega^{j-1} E_{j,p,q}=
\sum_{j=1}^{2q}  \omega^{j-1}
\log\frac{\Gamma\left(\frac{j+1}{2q}\right)}{\Gamma\left(\frac{j}{2q}\right)},$$
where 
\begin{equation}\label{even}
E_{j,p,q}  = \left\{\begin{array}{ll}
\log\frac{\mbox{\normalsize $\Gamma\left(\frac{j+1}{q}\right)$}}{\mbox{\normalsize $\Gamma\left(\frac{j}{q}\right)$}} 
& \mbox{if $p$ is even},\\
\log\frac{\mbox{\normalsize $\Gamma\left(\frac{j+1}{2q}\right)$}
\mbox{\normalsize $\Gamma\left(\frac{j+q}{2q}\right)$}}{\mbox{\normalsize $\Gamma\left(\frac{j}{2q}\right)$}
\mbox{\normalsize $\Gamma\left(\frac{j+q+1}{2q}\right)$}} & \mbox{if $p$ is odd}.
\end{array}\right.
\end{equation}
\end{theorem}
\noindent{\it Proof.}
By Definition~\ref{defff1} and formula (\ref{formc}), we only need to prove that 
\begin{equation}\label{ppo222}
-\sum_{n=1}^{\infty} \omega^{n-1}\log \frac{n+1}{n} =\sum_{j=1}^q \omega^{j-1} E_{j,p,q}
=\sum_{j=1}^{2q} \omega^{j-1}\log\frac{\Gamma\left(\frac{j+1}{2q}\right)}{\Gamma\left(\frac{j}{2q}\right)}.
\end{equation}
If $p$ is even,  $(\omega^{j-1})_{j \geq 1}$ 
is a periodic sequence of period $q$
with $\sum_{j=1}^q \omega^{j-1}=0$. In this case, 
the first equality in (\ref{ppo222}) follows immediately from Theorem~\ref{enana}. 

If $p$ is odd,  $(\omega^{j-1})_{j \geq 1}$ is a periodic sequence of period $2q$
with $\sum_{j=1}^{2q} \omega^{j-1}=0$. Hence in Theorem~\ref{enana} we may replace $q$ with $2q$, and
take $\theta_j=\omega^{j-1}.$
Since $\omega^{j+q}=-\omega^{j}$ for $j=1,2,\ldots,q$,
the first equality in (\ref{ppo222}) holds with $E_{j,p,q}$ given by (\ref{even}).

To show the second equality in (\ref{ppo222}), use Theorem~\ref{enana} and the fact that 
$\omega$ is a $2q\,$th root of unity. 
\hfill $\Box$

\vspace{.1in}

\begin{example}\label{five}
{\rm Take $p/q= 1/2$ and $2/3$. Using Theorem~\ref{identifu} and  identities 
(\ref{hajaja}) and (\ref{77sss}), we get
\begin{small}
\begin{eqnarray*}
\gamma(i) & = & \int_0^1\int_0^1 \frac{1-x}{(1-ixy)(-\log xy)}\,dx\,dy  
 =  \frac{\pi}{4}-\log\frac{\Gamma\left(1/4\right)^2}{\pi\sqrt{2\pi}}
+i \log\frac{8\sqrt{\pi}}{\Gamma\left(1/4\right)^2},\\
\gamma(e^{2\pi i/3})  & = &  \int_0^1\int_0^1 \frac{1-x}{(1-e^{2\pi i/3}xy)(-\log xy)}\,dx\,dy  = 
\frac{\pi}{4\sqrt{3}}-\frac{3}{2} \log\frac{2\pi}{3\Gamma\left(2/3\right)^2}
+i\left(\frac{\sqrt{3}}{2}\log\frac{9}{2\pi}-\frac{\pi}{12}\right).
\end{eqnarray*}
\end{small}
}
\end{example}


The following  corollary gives an additional formula for the value of $\gamma(z)$ 
at any point $\omega$ of the unit circle with argument $\pi p/q$, where $p$ and $q$ are
integers with $p$  odd and $q>0$. For such a point $\omega$, 
the formula expresses $\gamma(\omega)$ as a finite linear combination of powers of $\omega$ with real coefficients.

\begin{corollary}\label{extra}
Let $p$ and $q$ be integers, with $p$ odd and $q$ positive. If $\omega=e^{i\pi p/q}$, then
$$\gamma(\omega)=\sum_{j=1}^q \omega^{j-1} (D_{j,q}+E_{j,p,q}),$$
where $E_{j,p,q}$ is given by {\rm (\ref{even})}, and
\begin{equation}\label{22.5}
D_{j,q} =  \frac{1}{2q}\left[\psi\left(\frac{1}{2}+\frac{j}{2q}\right)-\psi\left(\frac{j}{2q}\right)\right].
\end{equation}
\end{corollary}
\noindent{\it Proof.} 
By Theorem~\ref{ss1111} and formula (\ref{formc}), we only need to show that
\begin{equation}\label{akks}
\sum_{n=1}^{\infty} \frac{\omega^{n-1}}{n}  =  \sum_{j=1}^q \omega^{j-1}D_{j,q}.
\end{equation}
Let $s_{1}, s_2, \ldots$ be the partial sums of the series.  Then, using $\omega^{2q}=1$,
$$\lim_{k \rightarrow \infty} s_{k}=\lim_{k \rightarrow \infty} s_{2qk}=
\lim_{k \rightarrow \infty} \sum_{j=1}^{2q} \omega^{j-1} \sum_{m=0}^{k-1}\frac{1}{2qm+j}.$$
Since $p$ is odd, $\omega^{j+q}=-\omega^{j}$ for $j=1,2,\ldots,q$. Therefore,
$$\lim_{k \rightarrow \infty} s_{k}=
\sum_{j=1}^q \omega^{j-1} \frac{1}{4q} \sum_{m=0}^{\infty} \frac{1}{\left(m+\frac{j}{2q}\right)\left(m+
\frac{j+q}{2q}\right)}.$$
Using the formula \cite[Section 8.363, Equation 3]{soso00}, \cite[Equation 44:5:8]{spanoer}
$$\sum_{m=0}^{\infty} \frac{1}{(m+x)(m+y)}=\frac{\psi(x)-\psi(y)}{x-y},$$
which is valid for $x,y \in \C-\{0,-1,-2,\ldots\}$ with $x\neq y$,
equation~(\ref{akks}) follows. \hfill $\Box$

\vspace{.1in}

\begin{remark}\label{oner}
{\rm In (\ref{22.5}) the quantity $\psi(j/q)$ can be calculated from Gauss's formula (\ref{gauss})
when $0 < j < q$, together with the value $\psi(1)=-\gamma$.}
\end{remark}


\begin{example}\label{six}
{\rm Setting  $p=1$ and $q=3$ in Corollary~\ref{extra}, we multiply by $6$ and obtain
\begin{eqnarray*}
6\gamma(e^{i\pi/3}) & = & \psi\left(\frac{2}{3}\right)-\psi\left(\frac{1}{6}\right)+6 \log 
\frac{\Gamma\left({1}/{3}\right)\Gamma\left({2}/{3}\right)}{\Gamma\left({1}/{6}\right)\Gamma\left({5}/{6}
\right)}  \\
& & +e^{i\pi/3}\left[\psi\left(\frac{5}{6}\right)-\psi\left(\frac{1}{3}\right)+6\log 
\frac{\Gamma\left({1}/{2}\right)\Gamma\left({5}/{6}\right)}{\Gamma\left({1}/{3}\right)}\right] \\
& & +e^{2i\pi/3}\left[-\gamma-\psi\left(\frac{1}{2}\right)+6\log 
\frac{6\Gamma\left({2}/{3}\right)}{\Gamma\left({1}/{2}\right)\Gamma\left({1}/{6}\right)}\right]\\
& = & \pi\sqrt{3}-3\log \frac{6\sqrt{3}}{\pi}+
i \left[\pi-3\sqrt{3} \log
\frac{\Gamma\left({1}/{6}\right)^2\Gamma\left({1}/{3}\right)^2\sqrt{3}}{24\pi^2}\right].
\end{eqnarray*}}
\end{example}

\vspace{.1in}

The next theorem gives a formula for the average value of the function $\gamma(z)$ at the vertices of
a regular polygon inscribed in the unit circle with the positive real axis a perpendicular 
bisector of one side.

\begin{theorem}\label{sumeee}
Given $q \in \Z^+$, let $\omega=e^{i\pi/q}$. Then
\begin{eqnarray*}
\frac{1}{q}\sum_{k=0}^{q-1} \gamma(\omega^{2k+1}) 
 & = & \int_0^1\int_0^1\frac{1-x}{(1+x^qy^q)(-\log xy )}\,dx \,dy  \\
& = & \frac{1}{2q}\left[\psi\left(\frac{q+1}{2q}\right)-\psi\left(\frac{1}{2q}\right)\right]-
\log\frac{\Gamma\left(\frac{1}{2q}\right)\Gamma\left(\frac{q+2}{2q}\right)}{\Gamma\left(\frac{1}{q}\right)\Gamma\left(
\frac{q+1}{2q}\right)}.
\end{eqnarray*}
\end{theorem}
We give two proofs.

\vspace{.1in}

\noindent {\it Proof 1.}
By the method of partial fractions, for $\zeta \in \C- \{\omega^{-1}, \omega^{-3}, \ldots,\omega^{-(2q-1)}\}$,
$$\sum_{k=0}^{q-1} \frac{1}{1-\zeta\omega^{2k+1}}=\frac{q}{1+\zeta^q}.$$
Letting $\zeta=xy$ and using (\ref{sonora1}), we obtain
\begin{equation}\label{444o}
\sum_{k=0}^{q-1} \gamma(\omega^{2k+1})=q \int_0^1\int_0^1\frac{1-x}{(1+x^qy^q)(-\log xy )}\,dx\,dy.
\end{equation}
Making the change of variables $X=x^{q}, Y=y^{q}$ gives
$$\sum_{k=0}^{q-1} \gamma(\omega^{2k+1})= \int_0^1\int_0^1\frac{(XY)^{\frac{1}{q}-1}-
X^{\frac{2}{q}-1}Y^{\frac{1}{q}-1}}{(1+XY)(-\log XY )}\, dX\,dY.$$
Now use the following evaluations  \cite[Corollary 3.4]{quisond}, which hold for $u,v>0$:
$$\int_0^1 \int_0^1 \frac{(XY)^{u-1}}{(1+XY)(-\log XY)}\,dX\,dY  =  
\frac{1}{2}\left[\psi \left(\frac{u+1}{2}\right)-\psi\left(\frac{u}{2}\right)\right],$$ 
$$\;\;\;\;\;\;\;\;\;\;\;\;\;\;\;\;\;\;\;\;\;\;\;
\int_0^1 \int_0^1 \frac{X^{u-1}Y^{v-1}}{(1+XY)(-\log XY)}\,dX\,dY  
=  \frac{1}{u-v}\log \frac{\Gamma\left(\frac{v}{2}\right)
\Gamma\left(\frac{u+1}{2}\right)}{\Gamma\left(\frac{u}{2}\right)\Gamma\left(\frac{v+1}{2}\right)}.
\;\;\;\;\;\;\;\;\;\;\;\;\;\;\;\;\;\;\;\;\;\;\mbox{\hfill $\Box$}$$

\vspace{.1in}

\noindent{\it Proof 2.} Using Corollary~\ref{extra} and the equality $E_{j,2k+1,q}=E_{j,1,q}$ from (\ref{even}), we obtain
$$\sum_{k=0}^{q-1} \gamma(\omega^{2k+1}) = \sum_{k=0}^{q-1}\sum_{j=1}^{q} (\omega^{2k+1})^{j-1}
(D_{j,q}+E_{j,2k+1,q})=
\sum_{j=1}^{q} (D_{j,q}+E_{j,1,q})\sum_{k=0}^{q-1} (\omega^{j-1})^{2k+1}.$$
Since $\omega=e^{i \pi/q}$, if $j \in \{1,\ldots,q\},$ then $\omega^{j-1} \neq -1$, and $\omega^{j-1}=1$ 
if and only if $j=1$.
It follows that $\sum_{k=0}^{q-1} (\omega^{j-1})^{2k+1}=0$ when $j\in\{2,\ldots,q\}$, hence
$$\sum_{k=0}^{q-1} \gamma(\omega^{2k+1}) =q(D_{1,q}+E_{1,1,q}).$$
The theorem now follows using (\ref{even}), (\ref{22.5}), and (\ref{444o}).  \hfill $\Box$

\vspace{.1in}

\begin{example}\label{seven}
{\rm Take $q=3$. Using (\ref{gauss}) and (\ref{77sss}), we get
$$\gamma(e^{i\pi/3})+\gamma(e^{i \pi})+\gamma(e^{5i\pi/3})=3 \int_0^1\int_0^1\frac{1-x}{(1+x^3y^3)(-\log xy )}\,dx\,dy =
\frac{\pi}{\sqrt{3}}-\log \frac{3\sqrt{3}}{2}.$$}
\end{example}

\vspace{.25in}

\noindent {\large {\bf 5. The Hyperfactorial $K$ Function and the Derivative of $\gamma(z)$}}\label{foursection}

\vspace{.25in}

In this section we study the derivative of $\gamma(z)$ at roots of unity, using
the following function.

\begin{definition}\label{3hyper}
{\rm The Kinkelin-Bendersky \cite{jinhhh,bend} {\it hyperfactorial $K$ function}  is the real-valued function $K(x)$ defined
for $x \geq 0$ by the relation
\begin{equation}\label{KKBB}
\log K(x)=\frac{x^2-x}{2}-\frac{x}{2}\log 2\pi+\int_0^x \log\Gamma(y)\,{\rm d}y.
\end{equation}}
\end{definition}
(Bendersky~\cite{bend} uses the notation $\Gamma_1$, but most authors today use $K$ to denote the function 
defined by~(\ref{KKBB}); see, for example,  \cite[p.\ 135]{FF}.) 

\vspace{.1in}

The $K$ function satisfies \cite[p. 279]{bend}
$K(0)=K(1)=K(2)=1$
and, for $x> 0$,
\begin{equation}\label{reveeee}
K(x+1)=x^xK(x).
\end{equation}
Thus, by induction on $n \in \Z^+$, we have
$K(n+1)=1^1 2^2 \cdots n^n.$

There is an analog for the $K$ function  \cite[p.\ 281]{bend} of Gauss's multiplication formula for the
gamma function. A special case is
$$\prod_{j=1}^{n-1} K\left(\frac{j}{n}\right)=\frac{A^{\frac{n^2-1}{n}}}{n^{\frac{1}{12n}}e^{\frac{n^2-1}{12n}}},$$
for $n\geq 2$, where  \cite[pp.\ 263-264]{bend}
\begin{equation}\label{norlal}
A=\lim_{x \rightarrow \infty} \frac{K(x+1)}{x^{\frac{x^2+x}{2}+\frac{1}{12}}\, e^{-\frac{x^2}{4}}}=
\lim_{n \rightarrow \infty} \frac{1^1 2^2 \cdots n^n}{n^{\frac{n^2+n}{2}+\frac{1}{12}}\, e^{-\frac{n^2}{4}}}=
1.28242712 \ldots
\end{equation}
is the {\it Glaisher-Kinkelin constant} (see \cite[p.\ 135]{FF}). In particular,
\begin{equation}\label{onehalf}
K\left(\frac{1}{2}\right)= \frac{A^{3/2}}{2^{1/24} e^{1/8}}.
\end{equation}

\vspace{.1in}

The following summation formula is needed in the proof of Theorem~\ref{0000ss1111}. 

\begin{theorem}\label{e224nana}
Let $(\theta_n)_{n \geq 1}$ be a periodic sequence in $\C$, and let $q$ be a positive multiple of its period.
If  \begin{equation} \label{thetajjj}
\sum_{j=1}^q \theta_j=0,
\end{equation}
then 
$$\sum_{n=1}^{\infty}   \theta_n \left(1-n\log\frac{n+1}{n}\right)=
\sum_{j=1}^q \theta_j \log \frac{K\left(\frac{j+1}{q}\right)^qe^{-j/q}}{K\left(\frac{j}{q}\right)^q\Gamma\left(\frac{j+1}{q}\right)}.$$
\end{theorem}
\noindent {\it Proof.}
Dirichlet's test implies that the series converges. 
The proof of the formula for its sum is similar to that of Theorem~\ref{enana}. The new ingredient is the 
{\it Barnes $G$ function}, which is defined  by the infinite product
\begin{equation}\label{43aa}
G(z+1)=\frac{(2\pi)^{z/2}}{e^{\frac{1}{2}z(z+1)+\frac{\gamma}{2} z^2}}\prod_{n=1}^{\infty} 
\frac{\left(1+\frac{z}{n}\right)^n}{e^{z-
\frac{1}{2n} z^2}}
\end{equation}
for $z \in \C$, and which is related to the $K$ 
function by $K(x+1)G(x+1)=\Gamma(x+1)^x$ for $x\geq 0$ (see \cite[Section 2.15]{FF}, 
\cite[Section 6.441]{soso00}). 

From (\ref{reveeee}) and (\ref{hajaja}), we can write the relation as $K(x)=G(x+1)^{-1}
\Gamma(x)^x$, so that
$$\frac{K\left(\frac{j+1}{q}\right)^q}{K\left(\frac{j}{q}\right)^q\Gamma\left(\frac{j+1}{q}\right)}=\left(
\frac{G\left(\frac{j}{q}+1\right)}{G\left(\frac{j+1}{q}+1\right)}\right)^q
\left(\frac{\Gamma\left(\frac{j+1}{q}\right)}{\Gamma\left(\frac{j}{q}\right)}\right)^j $$
for $j=1,2,\ldots,q$. Now take the logarithm, multiply by $\theta_j$, and sum from $j=1$ to $q$.
We compute the result in two parts. On the one hand, a calculation using (\ref{43aa}) and (\ref{thetajjj}) yields
$$\sum_{j=1}^q \theta_j q \log \frac{G\left(\frac{j}{q}+1\right)}{G\left(\frac{j+1}{q}+1\right)}=
\frac{1+\gamma}{q} \sum_{j=1}^q \theta_j j -\sum_{n=1}^{\infty} \sum_{j=1}^q \theta_j \left(\frac{j}{qn}+
qn\log \frac{qn+j+1}{qn+j}\right).$$
On the other hand, from (\ref{callclm}) we have
$$\sum_{j=1}^q \theta_j j \log \frac{\Gamma\left(\frac{j+1}{q}\right)}{\Gamma\left(\frac{j}{q}\right)}=
-\sum_{j=1}^q \theta_j j \left( \frac{\gamma}{q}+\log \frac{j+1}{j}\right)-\sum_{n=1}^{\infty} \sum_{j=1}^q 
\theta_j \left(-\frac{j}{qn}+j \log \frac{qn+j+1}{qn+j}\right).$$
Therefore, by addition
$$\sum_{j=1}^q \theta_j \log \frac{K\left(\frac{j+1}{q}\right)^q e^{-j/q}}{K\left(\frac{j}{q}\right)^q 
\Gamma\left(\frac{j+1}{q}\right)}=-\sum_{n=0}^{\infty} \sum_{j=1}^q \theta_j (qn+j)\log
\frac{qn+j+1}{qn+j}.$$
Substituting $\theta_j=\theta_{qn+j}$ on the right, and using (\ref{thetajjj}) again, we deduce the desired
formula. \hfill $\Box$

\vspace{.1in}

\begin{example}\label{eleven11}
{\rm Take $\theta_n=(-1)^{n-1}$ and $q=2$. Exponentiating the series, and 
using (\ref{reveeee}) and (\ref{onehalf}), we obtain the infinite product
$$\prod_{n=1}^{\infty} \left[\frac{e}{\left({1+\frac{1}{n}}\right)^{n}}\right]^{(-1)^{n-1}}=
\frac{2^{1/6}e\sqrt{\pi}}{A^6},$$}
\end{example}
where $A$ is the Glaisher-Kinkelin constant (\ref{norlal}).

\begin{remark}{\rm The last equality, in the form of the limits
$$\lim_{N\rightarrow \infty} \prod_{n=1}^{2N+1} \left(1+\frac{1}{n}\right)^{(-1)^{n+1}n}
=e \cdot \lim_{N\rightarrow \infty} \prod_{n=1}^{2N} \left(1+\frac{1}{n}\right)^{(-1)^{n+1}n}
=\frac{A^6}{2^{1/6}\sqrt{\pi}},$$
is due to Borwein and Dykshoorn~\cite{bowr}.}
\end{remark}

Our final theorem is the main result of the section. It gives the value of the derivative of 
$\gamma(z)$ at any point $\omega \neq 1$ of the unit circle with argument a rational multiple of~$\pi$.

\begin{theorem}\label{0000ss1111}
Let $p$ and $q$ be relatively prime positive integers. 
If $\omega=e^{i\pi p/q} \neq 1$, then the derivative of the function 
$\gamma(z)$ at $z=\omega$ is
$$\gamma'(\omega)=\frac{1}{\omega (1-\omega)}+\frac{\log(1-\omega)}{\omega^{2}}+\sum_{j=1}^{2q} 
\omega^{j-2} \log \frac{\Gamma\left(\frac{j}{2q}\right)
K\left(\frac{j+1}{2q}\right)^{2q}}{\Gamma\left(\frac{j+1}{2q}\right)^2K\left(\frac{j}{2q}\right)^{2q}}.$$
\end{theorem}
\noindent {\it Proof.}
Using Definition~\ref{defff1},
the Taylor series of the derivative of the product $z\gamma(z)$ when $|z|< 1$ is
\begin{equation}\label{kyparissovouno}
\gamma(z)+z\gamma'(z)=\sum_{n=1}^{\infty}z^{n-1}\left(1-n \log \frac{n+1}{n}\right).
\end{equation}
For any fixed $z \neq 1$ with $|z|=1$, the partial sums of the series $\sum_{n=1}^{\infty} z^{n-1}$ are bounded. It
follows, using Dirichlet's test, that the series in (\ref{kyparissovouno}) converges for such $z$.
Since the function $\gamma(z)+z\gamma'(z)$ is continuous at such $z$, Abel's limit theorem implies that
(\ref{kyparissovouno}) also holds when $|z|=1\neq z$.

In particular, we may take $z=\omega=e^{i \pi p/q}\neq 1$. Then
the sequence $(\omega^{n-1})_{n \geq 1}$ is periodic, and $2q$ is a positive multiple of its period. Since 
$\sum_{n=1}^{2q} \omega^{n-1}=0$, in Theorem~\ref{e224nana} we may replace $q$ with $2q$, and take 
$\theta_n=\omega^{n-1}$. 
Using Theorem~\ref{ss1111}, the result follows. \hfill $\Box$

\vspace{.1in}

Compare the following double integral formulas with those involving $A$
in \cite{quisond} and \cite{sonqqq}.

\begin{corollary} \label{gammaprime1}
The first and second derivatives of the function $\gamma(z)$  at $z=-1$ are
$$\gamma'(-1)=\int_0^1 \int_0^1 \frac{xy(1-x)}{(1+xy)^2(-\log xy)} \, dx\, dy=
\log\frac{2^{11/6}A^6}{\pi^{3/2} e}$$
and
$$\gamma''(-1)=\int_0^1 \int_0^1 \frac{2x^2y^2(1-x)}{(1+xy)^3(-\log xy)} \, dx\, dy=
\log\frac{2^{10/3}A^{24}}{\pi^{4} e^{13/4}}-\frac{7\zeta(3)}{2\pi^2}.$$
\end{corollary}
\noindent {\it Proof.} In both cases the first equality follows using Theorem~\ref{identifu}.

To obtain the value of $\gamma'(-1)$ use Theorem~\ref{0000ss1111} 
together with formulas~(\ref{reveeee}) and 
(\ref{onehalf}). Alternatively, use (\ref{kyparissovouno}), Example~\ref{eleven11}, 
and the value $\gamma(-1)=\log \frac{4}{\pi}.$ 

To evaluate $\gamma''(-1)$, differentiate the inversion formula in 
Theorem~\ref{basiliki} twice at $z=-1+it$, let $t$ tend to $0^+$, and re-arrange terms, to get 
$$\gamma''(-1)=4\gamma'(-1)-3\gamma(-1)+\frac{3}{4}-\gamma-\log \pi -\psi\left(\frac{1}{2}\right) +\frac{1}{4\pi^2}
\psi''\left(\frac{1}{2}\right)+i\left[\frac{1}{\pi} \psi'\left(\frac{1}{2}\right)-\frac{\pi}{2}\right].$$
Using the values \cite[Sections 44:8 and 44:12]{spanoer}
$$\psi\left(\frac{1}{2}\right)=-\gamma-\ln 4, \quad \psi'\left(\frac{1}{2}\right)=\frac{\pi^2}{2}, \quad 
\psi''\left(\frac{1}{2}\right)=-14\zeta(3),$$
the desired formula follows. \hfill $\Box$

\vspace{.1in}

\begin{center}
{\bf Acknowledgements}
\end{center}

The second author would like to thank Roger W. Barnard and Alexander Yu.\ Solynin for their advice 
in analysis.

\end{document}